\documentclass[letterpaper, 10 pt, conference]{ieeeconf}

\IEEEoverridecommandlockouts                              % This command is only needed if 
\usepackage{cite}
\usepackage{amsmath,amssymb,amsfonts}
\usepackage{algorithmic}
\usepackage{graphicx}
\usepackage{textcomp}

\usepackage{adjustbox}
\usepackage{todonotes}

\newtheorem{theorem}{\textbf{Theorem}}
\newtheorem{definition}{\textbf{Definition}}

\newtheorem{problem}{\textbf{Problem}}
\newtheorem{remark}{\textbf{Remark}}
\newtheorem{lemmaBold}{\textbf{Lemma}}
\newtheorem{example}{\textbf{Example}}

\title{\LARGE \bf
A Constrained Tracking Controller for Ramp and Sinusoidal Reference Signals using Robust Positive Invariance}

\author{Geovana Franca dos Santos$^{1}$ and Eugenio  B. Castelan$^{2}$ and Walter Lucia$^{1}$
\thanks{This work was supported in part by the Natural Sciences and Engineering Research Council of Canada (NSERC), {and in part by CNPq-Brazil, Grant 311567/2021-5.}}
\thanks{$^{1}$ Geovana Franca and Walter Lucia are with the CIISE Department, Concordia University, Montreal, QC, H3G-1M8, Canada {\tt\small geovana.francadossantos@mail.concordia.ca, walter.lucia@concordia.ca}.}
\thanks{$^{2}$ Eugênio Castelan is with the Department
of Automation and Systems Engineering, UFSC, Florianopolis, SC, 88040-900, Brazil  {\tt\small eugenio.castelan@ufsc.br}.}
}

\begin{document}

\maketitle
\thispagestyle{empty}
\pagestyle{empty}

% \todo[inline,color=yellow]{Title suggestion - Eugênio: Constrained tracking of ramp and sinusoidal reference signals via RPInvariance}

% \todo[inline, color=green]{Walter: I checked on this website "https://app.copyleaks.com/text-compare", what is the overlap with our L-CSS paper. I found that it is almost 60\%. Geovana go on that website, check which pieces are the same and try to edit/reword what is in the paper. With this high similarity score the paper might be rejected regardless of the content}

\begin{abstract}
This paper proposes an output feedback controller capable of ensuring steady-state offset-free tracking for ramp and sinusoidal reference signals while ensuring local stability and state and input constraints fulfillment. The proposed solution is derived by jointly exploiting the internal model principle, polyhedral robust positively invariant arguments, and the Extended Farkas' Lemma. In particular, by considering a generic class of output feedback controller equipped with a feedforward term, a proportional effect, and a double integrator, we offline design the controller's gains by means of a single bilinear optimization problem. A peculiar feature of the proposed design is that the sets of all the admissible reference signals and the plant's initial conditions are also offline determined. Simulation results are provided to testify to the effectiveness of the proposed tracking controller and its capability to deal with both state and input constraints.
\end{abstract}

% \begin{IEEEkeywords}
% Constrained tracking controller, bilinear programming, robust positive invariance.
% \end{IEEEkeywords}

\section{Introduction}
\label{sec:introduction}
The set invariance theory is widely used in control, particularly when dealing with constrained systems and stability analysis. Sets are the most appropriate language to specify several system performances, such as determining the domain of attraction or measuring the effect of persistent noise in a feedback loop. One of the most popular approaches in this field is the Lyapunov theory and positive invariance. A positively invariant set is a subset of the state space of a dynamical system with the property that, if the system state is in this set at some time, then it will stay in this set in the future \cite{blanchini2015set,tarbouriech2011stability}.

The Internal Model Principle (IMP)  \cite{chen2014linear} is an essential result for reference tracking problems. IMP provides the conditions under which a stabilizing controller also ensures tracking by assuming an unconstrained feedback control system. As a result, several constrained and unconstrained tracking controllers have been established in the literature from different perspectives. Present techniques vary from Proportional-Integral (PI) controllers to Model Predictive Control (MPC) \cite{carvalho2019multivariable}, as well as Reference and Command Governor solutions \cite{limon2005mpc,ferramosca2011optimal,di2015reference,garone2017reference}.

% \todo[inline, color=green]{Walter: the following paragraph (in red) is unclear to me:\\
% - It is unclear if \cite{tarbouriech2000output} uses ARE or LMIs.\\
% - are we categorizing by method or by plant's model? \\
% - why LMI appears twice? should it appear only once?} 

In particular, different PI-like control design methods have been developed for a variety of constrained systems that consider algorithms based on Linear Matrix Inequalities (LMIs). The controller developed in \cite {flores2008tracking} deals with linear time-invariant systems, the solution in \cite{figueiredoIFAC20} addresses Linear Parameter-Varying (LPV) systems, and the approach described in \cite{lopes2020anti} is designed for nonlinear systems represented by Takagi-Sugeno (TS) models. Also, linear systems with saturating actuators and additive disturbances using either Riccati equations (ARE) or LMIs are investigated in  \cite{tarbouriech2000output}, and uncertain linear systems under control saturation using a ``quasi"-LMI approach for periodic references have been considered in \cite{flores2009robust}. However, these solutions deal only with symmetrical input saturation constraints and define contractive ellipsoidal (or composite ellipsoidal) invariant areas. 
On the other hand, the PI-like controller with feedforward term proposed in \cite{geovana2023} leverages algebraic robust positive invariance relations to define a bilinear optimization design problem capable of simultaneously computing the controller's parameters and the set of admissible step reference signals. The proposed design methodology guarantees steady-state offset-free tracking, and it can deal with continuous-time linear systems subject to asymmetrical and polyhedral state and input constraints. 
%, which are  not feasible in \cite{tarbouriech2000output,flores2008tracking}. 
%Specifically, they use algebraic robust positive invariance relations to define a bilinear optimization design problem capable of simultaneously computing the PI controller parameters, the set of admissible reference signals, and the polyhedral state space region where the controller is robust positively invariant  \cite{blanchini2015set}. 
%In other words, the proposed controller ensures constraint fulfillment and boundedness of the state trajectory for any reference signal in the defined bound and for any initial condition inside the controller's domain. Moreover, it also guarantees steady-state offset-free tracking under any admissible constant reference signal. 
The discrete-time counterpart of the solution in \cite{geovana2023} is presented in \cite{sbai2023}.

 \subsection{Paper's contribution}
% \todo[inline, color=green]{Walter: \\
% - please check the text in blue which I have edited (the original text is commented below if you want to check/restore part of it).\\
% - I moved the discussion about the fact that we here consider only the single output case to a remark right after th model, see Remarl~\ref{remark:single-output-case} (do you agree/disagree?)\\
% - If you have in mind other contributions and/or differences compared to our previous works, please add them in this subsection}
% \color{blue}
In this manuscript, we extend the solution in \cite{geovana2023} (capable of dealing only with step reference signals) to design a controller that addresses a broader class of reference signals, including ramps and sinusoidal signals. The proposed controller's design procedure can deal with asymmetric and polyhedral state and input constraints. Differently from \cite{geovana2023}, we consider a generic second-order homogeneous representation of the exogenous reference signal and an enhanced PI-like controller structure with a double integrator term. Such a structure is then leveraged to define the augmented constrained dynamics of the controlled closed-loop system and design the controller using a single bilinear optimization program. The proposed solution has the peculiar feature of allowing the simultaneous design and optimization of the controller's parameters, controller's domain of attraction and set of admissible reference signals.

% Since, in practice, having a wider choice of reference and disturbance dynamics is often desirable, this study extends on the initial insights made in \cite{geovana2023} in continuous-time systems for tracking ramp and sinusoidal reference signals. Therefore, we propose a generic class of IMP-based output feedback controllers with additional proportional and feedforward terms to achieve this objective as an extension of the previous article. However, the main difference is the need for a double integrator instead of a single one to handle reference signals that are more complex than a set-point reference signal. As a result, following the same structure in \cite{geovana2023}, the order of the closed-loop system also increases by one dimension, and a bilinear optimization program allows the design of the controller's parameters using, for instance,  the nonlinear solver KNITRO \cite{knitro06}. 

%Numerical experiments and simulations will demonstrate the usefulness of this approach for monovariable ramp and sinusoidal references. Finally, for the sake of simplicity, this work analyses the scenario with a single input and a single output to track the reference signal. The proposed approach is also appropriate to consider the multivariable case, combining different reference signals among the input-output pairs. It will be the focus of future work and consider other signals of practical interest.
\color{black}

 The rest of the paper is organized as follows. Section \ref{sec:problem} proposes the constrained PI-like controller's design. Section \ref{sec:solution} presents the proposed solution, and Section \ref{sec:examples} shows a numerical example for ramp and sinusoidal reference signals. Finally, Section \ref{sec:conclusion} concludes the paper.

%%%%%%%%%%%%%%%%%%%%%%%%%%%%%%%%%%%%%%%%%%%%%%%%%%%%%%%%%%%%%%%%%%%%%%%%%%%%%%%%

\section{Constrained Control Problem}
\label{sec:problem}

Consider a Linear Time-Invariant Continuous-Time (LTIC) system in the form
\begin{equation}\label{eq:plant}
\begin{array}{rcl}
\dot{x}(t) & = & Ax(t)+B u(t), \\ 
 y(t) & = & C x(t),
 \end{array}
\end{equation}
where $x(t) \in \mathbb{R}^{n}$ is the state vector, $u(t) \in \mathbb{R}^{m}$ the control input vector, and $y(t) \in \mathbb{R}$ the measurement vector. The system matrices $(A, B, C)$ are of suitable dimensions, with $(A,B)$ controllable and $(C,A)$ observable.

% {\color{blue} the angular frequency $\omega$ of the sinusoidal reference signal does not assume complex conjugate poles}. 
%
%% pg 393 Blanchini's book
% \todo[inline]{Eugênio says :) Last assumption has to be verified! It may be valid for "ramps" but maybe not for sinusoidal signals !?!? }

\begin{remark}\label{remark:single-output-case}
For reference tracking purposes, in \eqref{eq:plant} we have assumed that the output of the system is scalar (i.e., $y(t)\in \mathbb{R}$). It is worth remarking that such a choice is not dictated by any intrinsic limitation of the solution hereafter proposed, but it is made only for the sake of clarity and to improve as much as possible the description and readability of the proposed solution.  
\end{remark}

\begin{definition} 
\label{def:RPI} 
\textit{A polyhedral set $\mathcal{P}(\phi) \subseteq \mathbb{R}^n$ is said to be Robust Positively Invariant (RPI) for the system $\dot{x}(t)=f(x(t),d(t)),$ $t\geq 0,$ $x(t) \in \mathbb{R}^n,$ $d(t) \in \Delta(\psi) \subseteq \mathbb{R}^{n_d}$,  where $\Delta(\psi)$ is a compact polyhedral set, if for any initial state $x(0)\in \mathcal{P}(\phi),$ the state trajectory $x(t)$ remains bounded inside $\mathcal{P}(\phi), \forall\, t \geq 0$ and $\forall\,  d(t) \in \Delta(\psi).$}
\end{definition}

The state and input vectors are assumed to be subject to the following state and input constraints 
\begin{equation}
\label{eq:x}
  x(t) \in \mathcal{X} =  \{x(t): Xx(t) \leq \textbf{1$_{l_x}$}\},~  X \in \mathbb{R}^{l_x \times n},
\end{equation}
\begin{equation}
\label{eq:u}
  u(t)\in \mathcal{U}  =  \{u(t): Uu(t) \leq \textbf{1$_{l_u}$}\},~ U \in \mathbb{R}^{l_u \times m}.
\end{equation}
%

% \todo[inline,color=green]{Walter: the proposed solution is for systems with just 1 output, right? Is this the meaning of monovariable that we discuss before? If yes, then somewhere we need to say that the proposed framework can deal also with system with more than one output and where each output can track a different reference signal. Nevertheless for clarity of exposition we here focus on the case where $p=1$}

For tracking purposes, we assume that $y(t)$ must track a  reference signal given by  the homogeneous equation
% \todo[inline]{Maybe we should associate a number to  this equation, in particular to refer when assuming that the system has also to be free of zeros, which are such that $\sigma^2 + \alpha \sigma = 0$.  }
%
\begin{equation}
    \ddot{r}(t)+\alpha r(t)  = 0,
\end{equation}

\noindent where the initial conditions are unknown and $r(t)$ describes 
% \todo[inline, color=green]{Walter: maybe (not sure, to be discussed), to be precise, we should also say the initial condition to ensure that we have a ramp with slop $a$ and a sinusoidal wave with amplitude ``a''? }
\begin{itemize}
     \item [\textit{i)}] 
     \textit{a ramp signal}, i.e., $r(t) = at,~t \geq 0,$ if $\alpha = 0$,
    \item [\textit{ii)}] 
     % \todo[inline, color=green]{Walter: We need just one of the conditions in red, right?}
     \textit{a sinusoidal signal}
    $r(t) =  a\sin \omega t,~t \geq 0,$ if $\alpha= \omega ^2$.
\end{itemize}

Also, the reference signal $r(t)$ is bounded in an asymmetric hyperrectangle described by the set
\begin{equation}
\label{eq:rt}
 \mathcal{R}(\rho)  =  \{r(t): Rr(t)\leq \rho\},
 \end{equation}
 
 \noindent where $ R = \begin{bmatrix} ~~I \\ -I\end{bmatrix} \in \mathbb{R}^{2\times 2}$ and $\rho = \begin{bmatrix} \rho_1  \\ \rho_2 \end{bmatrix} \in \mathbb{R}^{2}$.

Furthermore, the bounds of the ramp reference signal can be interpreted, for a given time $t$, as the slope $a=\dfrac{\rho}{t}$. For the sinusoidal signal, the bounds of $ \mathcal{R}(\rho)$ represent the amplitude $a = \rho$, although the angular frequency $\omega$ takes its place by appearing in the internal model's reference, which has to be part of the controller's structure. An additional condition that has to be assumed is $rank\begin{bmatrix} A-\sigma I_n & B\\ C & 0\end{bmatrix} = n+ 1$ for $\sigma$ a complex value, such that $\sigma^2 + \alpha  = 0$, i.e., the system is free from transmission zeros at the origin for $\alpha=0$ (ramp), or transmission zeros at $\pm j \omega $ for $\alpha = \omega^2$ (sinusoidal).
 
We assume that the tracking controller presents the following structure
\begin{equation}\label{eq:law}
u(t)  =  Ky(t)+K_{I_1}x_{I_1}(t)+K_{I_2}x_{I_2}(t)+K_rr(t),
\end{equation}
where $K,K_{I_1},K_{I_2}, K_r \in \mathbb{R}^{m}$, $e(t) = r(t) - y(t)$,
$$x_{I_1}(t) = \int_{0}^{t}(e(\tau)-\alpha x_{I_2})d\tau \in \mathbb{R}, $$
$$x_{I_2}(t) = \int_{0}^{t}x_{I_1}(s)ds \in \mathbb{R},$$
and $\alpha = 0$ for a ramp reference signal,  $\alpha = \omega ^2$ for a sinusoidal reference signal. 
%
% $$x_{I_1}(t) = \int_{0}^{t}e(\tau)d\tau \in \mathbb{R}, \text{ and }$$
% $$x_{I_2}(t) = \int_{0}^{t}\left(\int_{0}^{t}e(\tau)d\tau\right)ds \in \mathbb{R}.$$

\begin{remark}\label{remamrk:localIMP}
\it Note that $Ky$, $K_{I_1}x_{I_1}$ and $K_{I_2}x_{I_2}$ define a Proportional and double Integral effect, respectively, while $K_r$ is a feedforward term used to improve the reference response, see \cite[Chapter 5]{aastrom2006advanced}. 
Consequently,  according to the IMP \cite[Section 9.2.2]{chen2014linear}, any stabilizing controller having the structure of \eqref{eq:law}, guarantees asymptotic reference tracking. However, since the considered system is subject to state and input constraints, the validity of such a result may be restricted to a bounded state space region where the constraints are inactive.
\end{remark}

\begin{problem}\label{prob1}
\textit{Consider the constrained plant's model \eqref{eq:plant}-\eqref{eq:u}, the reference constraint \eqref{eq:rt} and the controller's structure \eqref{eq:law}. 
Design the control gains $(K,K_{I_1},K_{I_2}, K_r)$ in \eqref{eq:law}, the vector $\rho$ in \eqref{eq:rt}, and a RPI set $\mathcal{L}\subset \mathbb{R}^{n+2}$, such that for any initial condition  $x_{cl}(0) = [x^T(0) ~ x^T_{I_1}(0)  ~ x^T_{I_2}(0)]^T\in \mathcal{L},$ the reference $r(t)$ is asymptotically tracked and the constraints \eqref{eq:x}-\eqref{eq:u} are fulfilled.}
\end{problem}

A particular case of the previous problem considering piecewise set-point tracking references was investigated in \cite{geovana2023}. There, the control law \eqref{eq:law} was considered only with one integrator, that is, they do not have the gain $K_{I_2}$ and the integral state $x_{I_2}$. Then, the closed-loop system is the second-order instead of the third-order as in \eqref{eq:closed}.

% \todo[inline]{Eugênio says :) Maybe a remark (here or elsewhere) about our previous result(s) about set-point tracking, to emphasize the link and differences ?!?!}

%%%%%%%%%%%%%%%%%%%%%%%%%%%%%%%%%%%%%%%%%%%%%%%%%%%%%%%%%%%%%%%%%%%%%%%%%%%%%%%
\section{Proposed Solution}
\label{sec:solution}

Using the IMP and polyhedral robust positive invariance arguments, this section provides a solution for Problem \ref{prob1}. First, the closed-loop dynamics of \eqref{eq:plant} under the control input \eqref{eq:law} are considered, and the related constraints are explicitly stated. Then, using proper set inclusion conditions and the Extended Farkas' Lemma, all the necessary algebraic conditions characterizing the set of admissible controller's parameters, reference's bounds, and RPI sets $\mathcal{L}$ are derived. Finally, we described the resulting optimization problem for control design (see opt. \eqref{opt:nlin}).

The closed-loop system is described below
\begin{equation}
\label{eq:closed}
\dot{x}_{cl}(t)  =  A_{cl}{x_{cl}}(t) + B_{cl} 
 r(t),
\end{equation}
\noindent with $A_{cl}= \begin{bmatrix} A+BKC & BK_{I_1} & BK_{I_2}  \\ -C & 0 & -\alpha   \\ 0 & 1 & 0 \end{bmatrix}$, $B_{cl}=\begin{bmatrix} 
 BK_{r} \\  1 \\ 0  \end{bmatrix}$, $x_{cl}(t)=\begin{bmatrix} x(t) & x_{I_1}(t) & x_{I_2}(t) \end{bmatrix}^T$ and $n_{cl}=n+2$.

Note that the value of $\alpha$ in $A_{cl}$ changes depending on the reference signal.
%, in the case of the ramp $\alpha = 0$ and the sinusoidal signal $\alpha = \omega ^2$.

\begin{remark}
\textit{For the sake of notation clarity, in what follows, the dependency of $x,y,u,r$ from $t$ is omitted.}
\end{remark}

Furthermore, the input constraint \eqref{eq:u} is translated into a closed-loop constraint from \eqref{eq:law}, as follows
\begin{equation}
\label{eq:Uy}
 \mathcal{U}_{cl} = 
\left\{\begin{bmatrix} x_{cl} \\ r \end{bmatrix}: U \begin{bmatrix} KC & K_{I_1} & K_{I_2} & K_r \end{bmatrix} \begin{bmatrix} x_{cl} \\ r \end{bmatrix}\leq \textbf{1$_{l_u}$}\right\}.
\end{equation}

Since fast error-tracking dynamics are desirable, to minimize the magnitude of the vectors $x_{I_1}$ and $x_{I_2}$ impose a further optional constraint. In particular, we allow the possibility of bounding each component of $x_{I_1}$ and $x_{I_2}$ in the respective asymmetric interval
\begin{equation}
\label{eq:xi}
 \mathcal{X}_I= \left\{\begin{bmatrix}
    x_{I_1} \\ x_{I_2}
\end{bmatrix}: X_I\begin{bmatrix}
    x_{I_1} \\ x_{I_2}
\end{bmatrix} \leq \textbf{1}_{4}\right\}, 
\end{equation}
\noindent with $X_I=\begin{bmatrix}
X_{I_1} & 0   \\ 0 & X_{I_2}  \end{bmatrix} \in \mathbb{R}^{4 \times 2},~ X_{I_i} = \begin{bmatrix}
    ~~X_{I_{1i}} \\ -X_{I_{2i}}
\end{bmatrix}$, $i = 1,2$.
 
Thus, the set of state constraints acting on the closed-loop system (i.e., \eqref{eq:x} and \eqref{eq:xi}) can be re-written as the following single constraint
\begin{equation}
\label{eq:Xcl}
x_{cl}\in \mathcal{X}_{cl} = \left\{ x_{cl} :  X_{cl}  x_{cl}  \leq \textbf{1$_{l_{x_{cl}}}$}\right\},
\end{equation}
with $X_{cl}=\begin{bmatrix}
X & 0 \\ 0 & X_I
\end{bmatrix} \in \mathbb{R}^{l_{x_{cl}} \times n_{cl}}$, $l_{x_{cl}}=l_{x}+l_{x_{I_1}}+l_{x_{I_2}}.$

Given that the IMP is only locally valid for the constrained system, the idea is to characterize an RPI polyhedral set $\mathcal{L}$ for \eqref{eq:closed}, where the state trajectory $x_{cl}$ is confined and constraints \eqref{eq:Uy} and \eqref{eq:Xcl} are fulfilled for any admissible reference signal.

The following comment is similar to Remark 4 in \cite{geovana2023}, adapted for the closed-loop system and the reference signals considered in the present work.

\begin{remark}
\it 
The reference signal $r$ in \eqref{eq:closed} can be interpreted as a bounded disturbance and, consequently, the RPI nature of $\mathcal{L}$ can be described from Definition~\ref{def:RPI} substituting 
$\dot{x} = f(x,d) \leftarrow \eqref{eq:closed}$, $\mathcal{P}(\phi) \leftarrow \mathcal{L},~ d \leftarrow r, \text{ and } \Delta(\psi) \leftarrow \mathcal{R}(\rho).
$
\end{remark}

Let $\mathcal{L}$ be described by the polyhedral set
\begin{equation}
\label{eq:setL}
\mathcal{L} = \{x_{cl}:L_{cl}x_{cl}\leq \mathbf{1}_l\},
 \end{equation}
with $L_{cl}=\begin{bmatrix}
    L & L_{I_1} & L_{I_2}
\end{bmatrix} \in \mathbb{R}^{l \times n_{cl}}$ and $rank(L_{cl})=n_{cl},$  it is possible to state  Theorem~\ref{prop:solution} that defines the algebraic conditions under which the controller \eqref{eq:law} provides a solution to Problem~\ref{prob1}. 

% \todo[inline]{Suggestion: present the Th. as a Lemma, and maybe change Proposition to Theorem. Put in the appendix the Lemmas's proof and bring back the Th/Proposition proof to the text.}

\begin{lemmaBold}
\label{lemma:RPI}
A polyhedron $\mathcal{L}$, \eqref{eq:setL}, is a RPI set of the system \eqref{eq:closed}, if and only if, exists a Meztler matrix $H \in \mathbb{R}^{l \times l}$, the matrix $H_r>0 \in \mathbb{R}^{l \times l_r}$ and a scalar $\gamma>0$, such that 
\begin{equation}
\label{eq:deltainv}
\begin{array}{rcl}
HL & = & L(A+BKC)-L_{I_1}C, \\
HL_{I_1} & = & LBK_{I_1}+L_{I_2}, \\
HL_{I_2} & = & LBK_{I_2}-L_{I_1}\alpha, \\
H_{r}R & = & LBK_r+L_{I_1},\\
H\textbf1_{l}+H_{r}\rho &\leq& - \gamma 1_{l}.
\end{array}
\end{equation}
\end{lemmaBold}

Remember that the value of $\alpha$ is defined based on the reference signal, being null for the ramp and equal to $\omega ^2$ for the sinusoidal signal. Also, it is important to point out that the matrix $L_{cl}$ in \eqref{eq:setL} has a complete column rank, if and only if, it admits pseudo-inverse matrices  $V_1 \in \mathbb{R}^{n \times l}$, $V_2, V_3 \in \mathbb{R}^{1 \times l}$ so that
\begin{equation}
\label{eq:rank}
\begin{bmatrix}
    V_1 \\ V_2 \\ V_3
\end{bmatrix} \begin{bmatrix}
    L & L_{I_1} & L_{I_2}
\end{bmatrix} = I_{n_{cl}}.
\end{equation}

The necessary and sufficient algebraic conditions to obtain the inclusions between the polyhedral sets involved $\mathcal{X}_{cl}$ and $\mathcal{U}_{cl}$, respectively, can be obtained by applying an extension of Farkas' Lemma presented in Appendix A, as follows:

\begin{itemize}
    \item Assume that there exists a matrix $L_{cl} \in \mathbb{R}^{l \times n_{cl}},$ with $l > n_{cl}$,  non-negative matrices $T_1 \in \mathbb{R}^{l_x \times l}$, $T_2 \in \mathbb{R}^{l_{x_{I_1}} \times l}$, $T_3 \in \mathbb{R}^{l_{x_{I_2}} \times l}$,
and a scalar $\gamma>0$ satisfying
\begin{equation}
\begin{array}{rcl}
\label{eq:state}
% TL_{cl}
\begin{bmatrix}
    T_1 \\ T_2 \\ T_3
\end{bmatrix} \begin{bmatrix}
    L & L_{I_1} & L_{I_2}
\end{bmatrix}
&= & \begin{bmatrix}
X & 0 & 0 \\
0& X_{I_1} & 0   \\ 0& 0 & X_{I_2}  \end{bmatrix} ,  \\
T_1 \textbf{1}_l & \leq & \textbf{1}_{l_{x}} \\
 T_2 \textbf{1}_l & \leq & \textbf{1}_{l_{x_{I_1}}}, \\
 T_3 \textbf{1}_l & \leq & \textbf{1}_{l_{x_{I_2}}}.
\end{array}
\end{equation}

\item There exists a non-negative matrix $ Q \in \mathbb{R}^{l_u \times l}$ and $ Q_{r} \in \mathbb{R}^{l_u \times l_r}$, such that
   \begin{equation}
\begin{array}{rcl}
\label{eq:control}
QL & = & UKC,  \\
QL_{I_1} & = & UK_{I_1}, \\
QL_{I_2} & = & UK_{I_2},  \\
Q_{r}R & = & UK_{r}, \\
Q \textbf{1}_l + Q_{r}\rho & \leq & \textbf{1}_{l_u}.
\end{array}
\end{equation}
\end{itemize}

\begin{theorem}\label{prop:solution} \textit{Consider the closed-loop system \eqref{eq:closed}, the polyhedral sets \eqref{eq:rt}, \eqref{eq:Uy} and \eqref{eq:Xcl}. Assume that the Lemma \ref{lemma:RPI} and the condition \eqref{eq:rank}-\eqref{eq:control} are satisfied. Then, the polyhedral set $\mathcal L$ is RPI, such that $\mathcal{L} \subseteq \mathcal{X}_{cl}$ and  $\begin{bmatrix}KC &  K_{I_1} & K_{I_2} \end{bmatrix}  \mathcal{L} \oplus K_r \mathcal{R}(\rho)\subseteq \mathcal{U}_{cl}$, where $\oplus$ denotes the Minkowski set sum operator.} Therefore, for any initial condition $x_{cl}(0) = [x^T(0) ~ x^T_{I_1}(0)  ~ x^T_{I_2}(0)]^T \in \mathcal L$ the output $y$ asymptotically tracks any reference $r \in \mathcal R (\rho)$, with corresponding closed-loop trajectories fulfilling  the prescribed constraints. 
\end{theorem}

 {\bf \textit{Proof:}} The proof is adapted from  \cite{geovana2023} as follows:

 First, condition \eqref{eq:rank} is equivalent to imposing  $rank(L_{cl}) = n_{cl}$ and, possibly, that $\mathcal L$ is compact. Then, to guarantee that the system will stay within the closed loop state constraints, we impose the inclusion $\mathcal{L} \subseteq \mathcal{X}_{cl}$, which, by applying the Extended Farkas' Lemma \ref{lemma:EFL}, is equivalent to the existence of the non-negative matrix $T = \begin{bmatrix}
    T_1^T & T_2^T & T_3^T
\end{bmatrix}^T$, that satisfies the relation \eqref{eq:state}. Likewise, applying Lemma \ref{lemma:EFL}, the existence of non-negative matrices $Q$ and $Q_r$ verifying \eqref{eq:control}, represents the inclusion $\begin{bmatrix}KC & K_{I_1} & K_{I_2} \end{bmatrix}  \mathcal{L} \oplus K_r \mathcal{R}(\rho)\subseteq \mathcal{U}$ or, equivalently, 
\begin{equation*}
\begin{array}{rcl}
\left[\begin{array}{cc}Q & Q_r \end{array}\right]\left[\begin{array}{cc}L_{cl} & 0 \\ 0 & R \end{array}\right]& = & U\left[\begin{array}{cccc} KC & K_{I_1} & K_{I_2} & K_{r}\end{array}\right], \\
\begin{bmatrix}Q & Q_r \end{bmatrix}\begin{bmatrix}\textbf{1}_l \\\rho \end{bmatrix} & \leq & \textbf{1}_{l_u}.
\end{array}\end{equation*}
Finally, under the conditions  \eqref{eq:deltainv}-\eqref{eq:control}, for any reference signal $r\in \mathcal R (\rho)$ and for all $x_{cl}(0) = [x^T(0) ~ x^T_{I_1}(0)  ~ x^T_{I_2}(0)]^T \in \mathcal L,$ the closed-loop state trajectory remains inside $\mathcal{L}$ while fulfilling all the prescribed state and input constraints. Consequently, the system evolves in a domain where all the constraints are inactive and the closed-loop dynamics are uniquely determined by the unconstrained linear model \eqref{eq:closed}, whose state matrix $A_{cl}$ is Schur stable. Hence, the IMP is locally valid for any ramp or sinusoidal reference signal $r \in \mathcal R (\rho)$ and for all $x_{cl}(0) = [x^T(0) ~ x^T_{I_1}(0)  ~ x^T_{I_2}(0)]^T \in \mathcal L$.  \hfill $\Box$

The particular and simplest case for set-point reference tracking problem was covered in \cite{geovana2023} and compared with the LMI-based solution developed in \cite[Sec. 8.6.2]{blanchini2015set}. Note that Theorem 1 considers the monovariable case for clarity
 as a preliminary result, which does not prevent the application of this approach to the multivariable case, for example, the combination of a step, ramp, or sinusoidal reference signal, which will be considered in future work.

% \todo[inline]{Eugênio says :) I think that it is necessary to give some remarks about the consequences of the Proposition to each case, ramp and sinusoidal reference. Otherwise, it can be viewed as "you are doing the same as you did before"! For instance, the interpretation of the bounds for the ramp involves the $at$, but for sinusoidal signals, it seems be mostly related to the amplitude, although the frequency takes it place  by appearing in the internal model. Other "questions": - Can we consider piece-wise ramps (kind of triangular shapes?!!) or sinusoidal references? - The ramp problem also admits step signals, so the combination??? Some of this questions may be left as open questions or possibilities, but mentioning them can "make the difference(s)"!}

Next, let us consider the set of decision variables for the design of the controller \eqref{eq:law} given by
\begin{equation*}\label{eq:decision_variables}
    \lambda(\cdot) = (K,K_{I_1},K_{I_2},K_r,L_{cl}, H,H_{r},T, Q,  Q_{r}, V, X_I, \lambda,\rho).
\end{equation*}

The algebraic relations \eqref{eq:deltainv}-\eqref{eq:control} define the constraints under which the controller provides a solution to Problem~\ref{prob1}. The resulting bilinear optimization problem is presented below.
\begin{equation}
\label{opt:nlin}
\begin{aligned}
& \underset{\lambda(\cdot)}{\text{maximize}}
& & \Phi(\cdot),\\
&  \text{subject to}
& &  \eqref{eq:deltainv} - \eqref{eq:control},  \\
& & &f_\ell(\cdot) \leq \varphi_\ell, ~\ell=1,\ldots,\bar{\ell},\\
\end{aligned}
\end{equation}
where $\Phi(\cdot)$ is the cost function and  $f_\ell(\cdot)\leq \varphi_\ell$ are $\bar{\ell}$ auxiliary constraints instrumental to imposing limits over all the non-bounded decision variables. 

Furthermore, the choice of the cost function $\Phi(\cdot)$ depends on the designer's objectives. There are two options:
\begin{itemize}
    \item [\textit{i)}] $\Phi(\cdot)=\Phi_1 = \rho_1+\rho_2$, which allows us to maximize the hyperrectangle $\mathcal{R}(\rho)$ of all admissible reference signals.  
    \item [\textit{ii)}] $\Phi(\cdot)= \Phi_2 = - trace(X_{I_{11}}+X_{I_{12}}+X_{I_{21}}+X_{I_{22}})$, which allows us
    to minimize the limits of the admissible integral errors $\mathcal{X}_I$.
\end{itemize}

Notice that the opt.~\eqref{opt:nlin} is bilinear because it involves multiplication between decision variables, and therefore, it can be solved by nonlinear optimization techniques. One possible way to solve the bilinear optimization \eqref{opt:nlin} is to resort to the nonlinear state-of-the-art solver KNITRO \cite{knitro06}, which has already been successfully employed for similar problems in  \cite{dos2021pi,briao2021output,ernesto2021incremental,martins2020metodo}.  Thus, one can find the bounds of the decision space using insights about the plant's constraint limits and a trial-and-error approach. For further discussions about KNITRO and its use, see \cite[Section 4.2]{briao2021output}.

Finally, we summarize below the number of variables and constraints characterizing~\eqref{opt:nlin}. The proposed solution's complexity increases with the plant's dimensions, state and input constraints, reference set, and RPI set $\mathcal{L}$ complexity  \cite{geovana2023}:
\begin{itemize}
    \item \textbf{\# of variables} = $m+l(n_{cl}+l+l_r+l_{x}+l_{x_{I_1}}+l_{x_{I_2}}+l_u)+2(l_u+2)+n_{cl}^2+1$,
    \item \textbf{\# of equalities} =  $n_{cl}(l+l_{x}+l_{x_{I_1}}+l_{x_{I_2}}+l_u+n_{cl})+2(l+l_u)$,
    \item \textbf{\# of inequalities} = $l+ l_{x}+l_{x_{I_1}}+l_{x_{I_2}}+l_u$.    
\end{itemize}

%%%%%%%%%%%%%%%%%%%%%%%%%%%%%%%%%%%%%%%%%%%%%%%%%%%%%%%%%%%%%%%%%%%%%%%%%%%%%%
\section{Simulation}
\label{sec:examples}

In this section, the proposed tracking controller design's effectiveness is validated through simulation results obtained for a linearized model of a two-tank system. In the first case, a piecewise ramp signal is the reference input, and the controller's performance is evaluated for the two proposed cost functions. For the second case, we consider the sinusoidal signal as input, and for the sake of space, we showed the results only for the second cost function.

In order to reduce the search space and improve the numerical performance using the nonlinear solver KNITRO \cite{knitro06}, the optimization variables in \eqref{opt:nlin} have been bounded (element by element) as $H, H_r, T, Q, Q_r, G \text{ in } \begin{bmatrix}0, 10^2\end{bmatrix}$, $L,K, K_{I_1}, K_{I_2}, K_r \text{ in } \begin{bmatrix}-10^2, 10^2\end{bmatrix}$, and $V \text{ in } \begin{bmatrix}-10^3, 10^3\end{bmatrix}$.

\begin{example}\label{exammple:tank} Consider the following linearized model of a two-tank system (adapted from \cite{ferramosca2011optimal}):
\end{example}
\begin{equation}
\label{eq:ex1}
\begin{array}{rcl}
\dot{x} & = & \begin{bmatrix}
-0.0304 &~~ 0.0187 \\ 0 & -0.0187
\end{bmatrix} x+\begin{bmatrix}
6.6667 \\ 10
\end{bmatrix} u, \\  
 y & = & \begin{bmatrix}
 1 & 0 
 \end{bmatrix} x,
 \end{array}
\end{equation}
where the state and input constraints are: $-0.38 \leq x_1 \leq 0.68$, $-0.35 \leq x_2 \leq 0.65$, and $\|u\|_\infty \leq 2.$ 
We have solved the optimization problem \eqref{opt:nlin} for both $\Phi(\cdot)=\Phi_1$ and $\Phi(\cdot)=\Phi_2$ using $l=9$ for the following piecewise ramp signal:
\begin{equation}
\label{eq:ref_ramp}
r = \left\{ \begin{array}{lcl} 
~~0.01t &\text{IF}& t \leq 30\,\sec,\\
-0.0071t+0.5143 &\text{IF}& 30<t \leq 100\,\sec, \\
-0.2 &\text{IF}& t > 100\,\sec.\end{array}\right.
\end{equation}

Table \ref{table:comparison_example1} summarizes and compares the design results for both objective functions. In particular, $\rho$ defines the bounds of the set of admissible reference signals $\mathcal{R}(\rho)$, $X_I$ is the shaping matrix of the set constraining the integral error, the gains $K,K_{I_1},K_{I_2}$ and $K_r$ define the control law \eqref{eq:law}. As expected, for the cost function $\Phi_1,$ the bounds for the admissible reference signals $\mathcal{R}(\rho)$ are bigger than the ones obtained for $\Phi_2$. On the other hand, by using $\Phi_2,$ it is possible to obtain a smaller integral error and faster reference tracking. Indeed, for $\Phi_1,$ $x_{I_1} \in [-14 ,\, 16],$ and $x_{I_2} \in [-21 ,\, 33],$ while for $\Phi_2,$ $x_{I_1} \in [-10,\, 10]$  and $x_{I_2} \in [-13 ,\, 20]$. However, the cost to pay is a reduced size for the set of admissible reference $\mathcal{R}(\rho).$ 

\begin{table}[ht!]
\centering
\caption{Design results using opt.~\eqref{opt:nlin}: $\Phi_1$ vs $\Phi_2$}\label{table:comparison_example1}
\begin{tabular}{|c|c|c|c|c|}
\hline
$\Phi_i$ & $\rho$ & $X_I$ & $[K~ K_{I_1} ~ K_{I_2}~ K_r]^T$  \\ \hline
1        & $\begin{bmatrix} 0.4682 \\ 0.1606\end{bmatrix}$  & $\begin{bmatrix} ~~0.0617 \\   -0.0693 \\ ~~0.0300 \\ -0.0456\end{bmatrix}$ & $\begin{bmatrix}-3.3170 \\  0.3141 \\  0.0071 \\  2.8208\end{bmatrix}$   \\ \hline
2        & $\begin{bmatrix}0.3000 \\ 0.2000 \end{bmatrix}$  & $\begin{bmatrix} ~~0.1000 \\ -0.1000 \\ ~~0.0487 \\ -0.0739 \end{bmatrix}$   & $\begin{bmatrix} -3.8881 \\   0.3733 \\ 0.0085 \\ 3.3142\end{bmatrix}$  \\ \hline
\end{tabular}
\end{table}

\begin{figure}[ht!]
    \centering  
      \includegraphics[width=0.85\columnwidth]{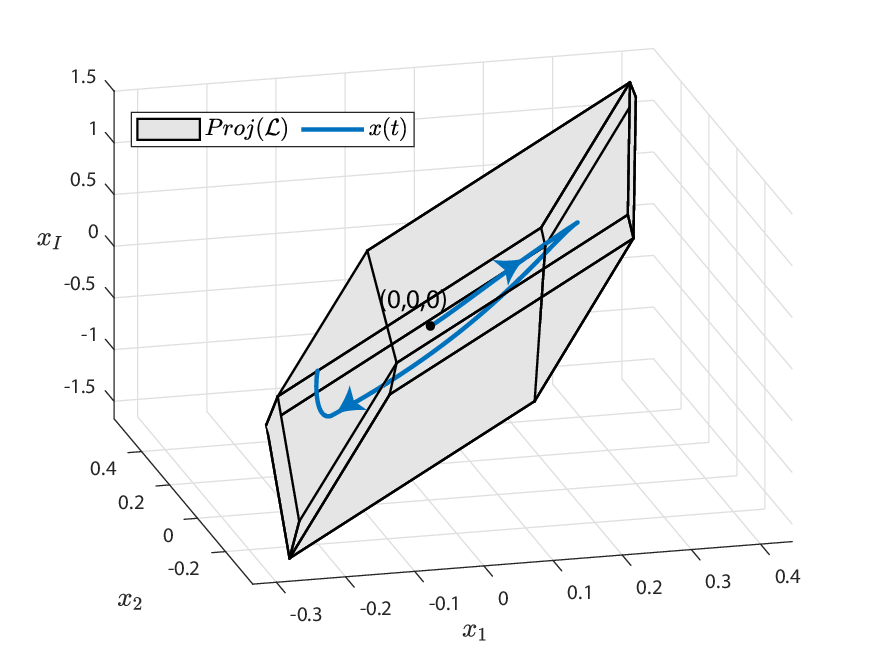} 
        \caption{Projected RPI set $\mathcal{L}$ and state trajectory for reference \eqref{eq:ref_ramp}}
        \label{fig:ramp}
\end{figure}

\begin{figure}[ht!]
    \centering       
    \includegraphics[width=0.8\columnwidth]{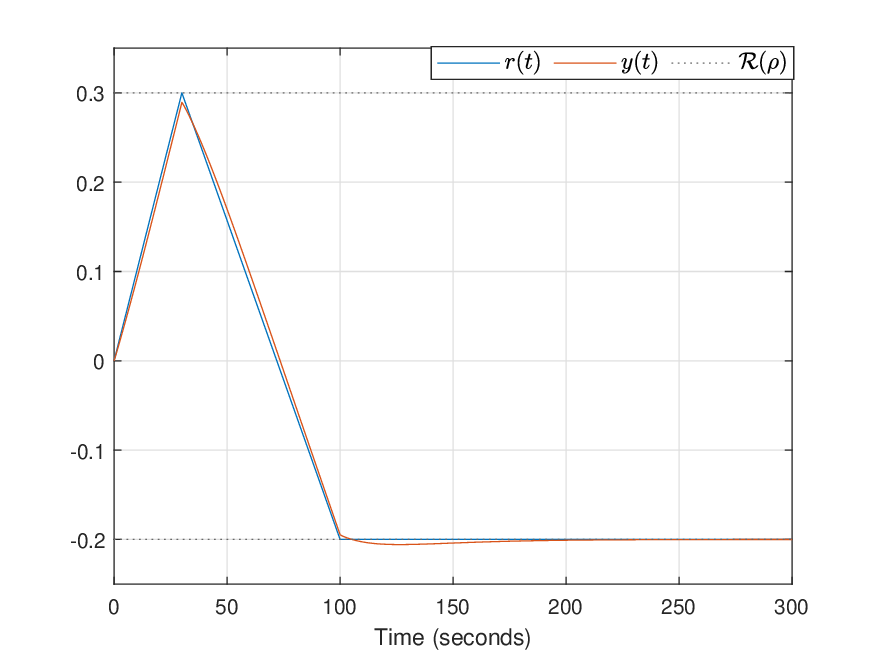} 
        \caption{Output and reference tracking for reference \eqref{eq:ref_ramp}}
        \label{fig:ramp_resp}
\end{figure}

\begin{figure}[ht!]
    \centering  
     \includegraphics[width=0.8\columnwidth]{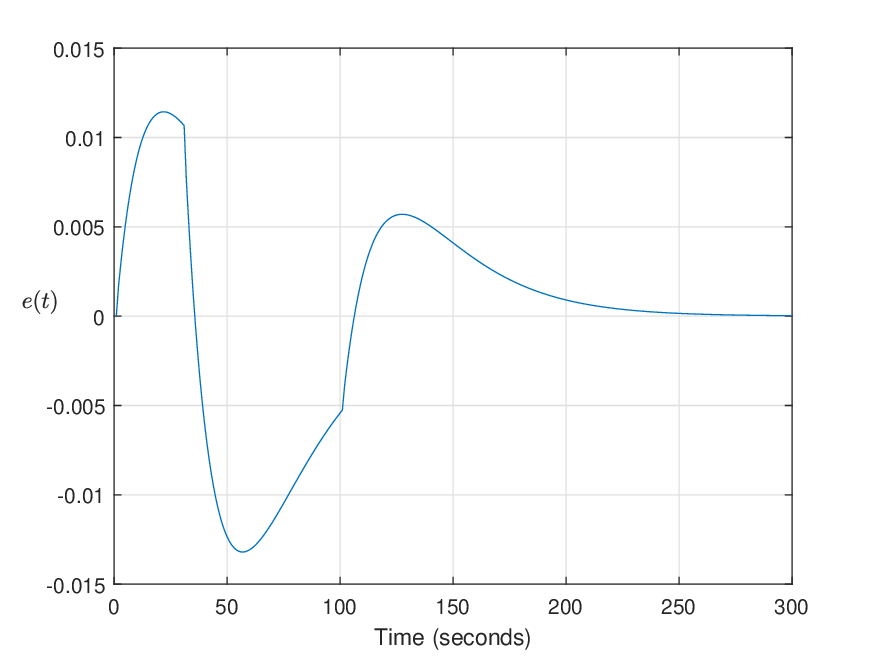} 
        \caption{Error signal for reference \eqref{eq:ref_ramp}}
        \label{fig:ramp_control}
\end{figure}

The projection of the resulting RPI set $\mathcal{L}$ is depicted in Fig. \ref{fig:ramp} with the closed-loop trajectorys obtained using the tracking controller associated to $\Phi(\cdot)=\Phi_2.$ It has been obtained starting from a zero initial condition and following the piecewise ramp trajectory $r$ in \eqref{eq:ref_ramp}, as shown in Fig. \ref{fig:ramp_resp}. The obtained trajectory confirms that the designed tracking controller allows the plant to asymptotically track the assigned piecewise ramp signal while ensuring that the state trajectory remains confined in the constraint-admissible RPI set $\mathcal{L}.$ Fig. \ref{fig:ramp_control} shows the error signal for the piecewise ramp reference. Note that, as expected, the error approaches zero when the plant reaches a steady-state regime, i.e., after the the reference signal is unchanged for a sufficiently long time period.

\begin{figure}[ht!]
    \centering  
      \includegraphics[width=0.85\columnwidth]{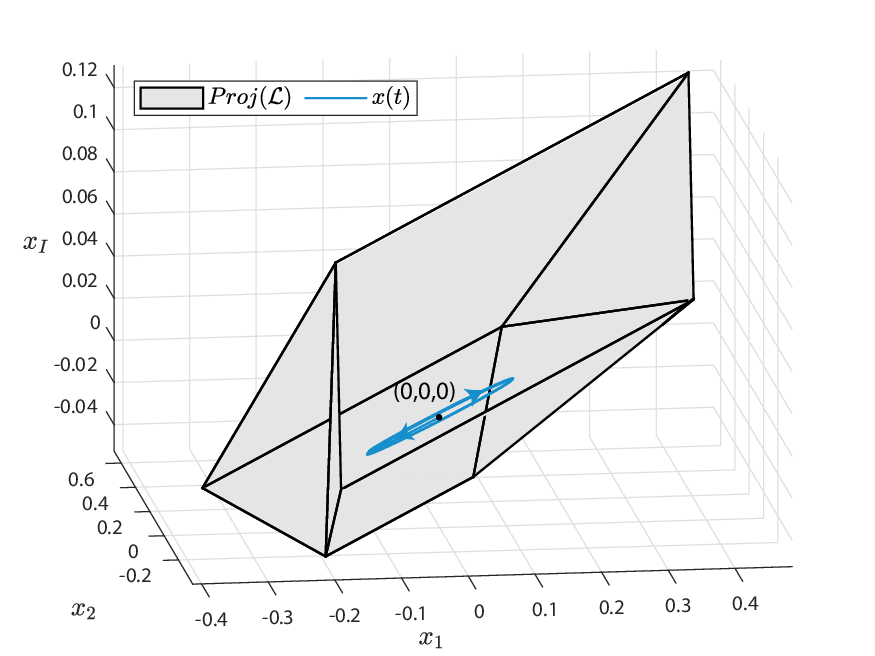} 
        \caption{Projected RPI set $\mathcal{L}$ and state trajectory for sinusoidal reference}
        \label{fig:sine}
\end{figure}

\begin{figure}[ht!]
    \centering  
      \includegraphics[width=0.8\columnwidth]{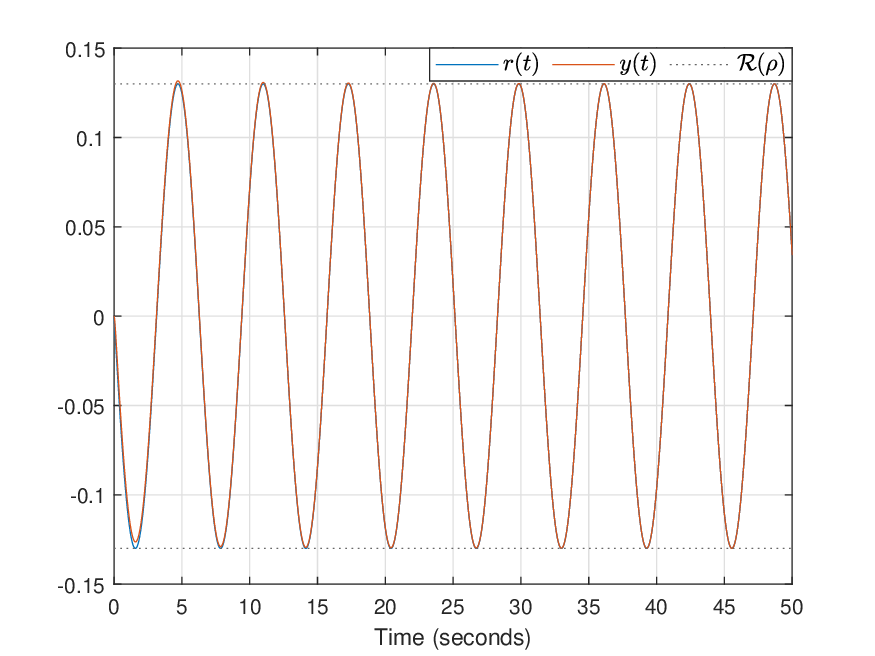} 
        \caption{Output and reference tracking for sinusoidal reference}
        \label{fig:sine_resp}
\end{figure}

\begin{figure}[ht!]
    \centering  
     \includegraphics[width=0.8\columnwidth]{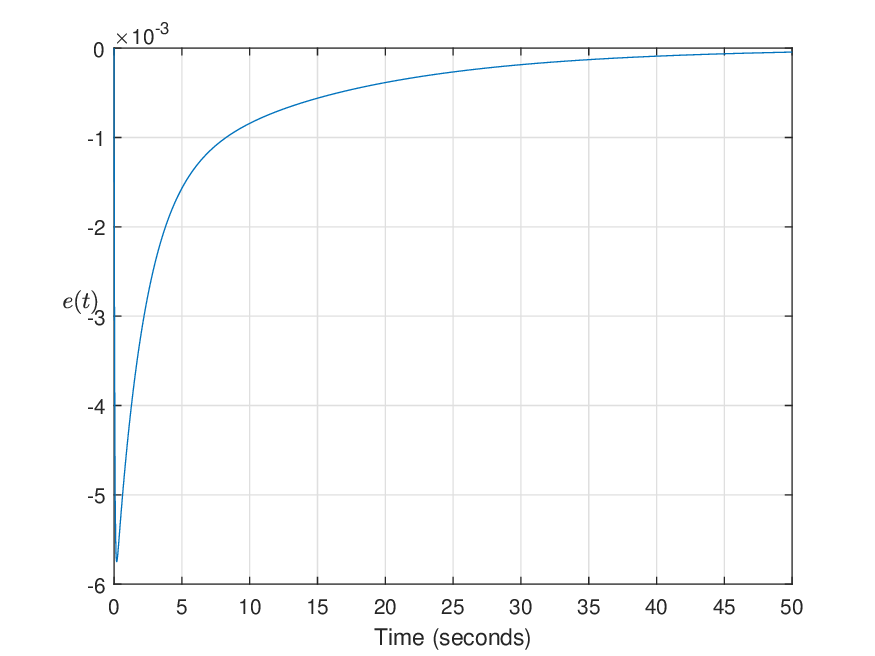} 
        \caption{Error signal for sinusoidal reference}
        \label{fig:sine_control}
\end{figure}

Next, we define a sinusoidal signal as a reference input $r =  a\sin \omega t$, considering $\omega =1$ and the objective function $\Phi(\cdot)=\Phi_2$ in the optimization problem \eqref{opt:nlin}. By solving \eqref{opt:nlin}, the following results have been obtained: $K = -3.2391$, $K_{I_1} =  1.6503$, $K_{I_2} =  -3.1178$, $K_r = 3.2106$, $X_I= \begin{bmatrix}0.1000 & -0.1000 & 0.1000 & -0.1000 \end{bmatrix}^T$ and $a = \rho = 0.13$. Fig. \ref{fig:sine} shows the projection of the RPI set $\mathcal{L}$ and the trajectory starting from a zero initial condition and following the sinusoidal signal, as shown in Fig. \ref{fig:sine_resp}, respecting the reference constraints $\mathcal{R}(\rho)$. The error signal  shown in Fig. \ref{fig:sine_control} confirms that the proposed controller ensures steady state offset-free tracking of the considered sinusoidal reference signal. 

%%%%%%%%%%%%%%%%%%%%%%%%%%%%%%%%%%%%%%%%%%%%%%%%%%%%%%%%%%%%%%%%%%%%%%%%%%%%%%%%
\section{Conclusion}
\label{sec:conclusion}

In this paper, we have extended and generalized the approach in \cite{geovana2023} to design tracking controllers for constrained linear time-invariant systems whose output is required to track ramp and sinusoidal reference signals. By assuming polyhedral state and input constraints, we have jointly resorted to robust positively invariant sets theory and IMP to define a single bilinear programming problem, which solves the control design problem. The proposed solution can simultaneously compute the controller parameters, the set of admissible reference signals, and the controller's domain of attraction. 

Future works will be devoted to extending the proposed approach to deal with other classes of reference signals and with the possibility that each system output (in multi-output setups) must track a differently-shaped reference trajectory.

%Our proposed extension allows tracking ramp or sinusoidal reference signals but with limits imposed by the state and control constraints. In the future, we plan to generalize this approach to multivariable systems. This will also associate different classes of reference signals to the input-output pairs.

\section*{Appendix A}

This appendix recalls basic definitions for polyhedral sets and the Extended Farkas' Lemma adapted from \cite{hennet1995discrete,blanchini2015set}.

\begin{definition} 
\label{def:polyhedral}
\textit{Any closed and convex polyhedral set $\mathcal{P}(\phi) \subseteq \mathbb{R}^n$ can be characterized by a shaping matrix $P \in \mathbb{R}^{l_p \times n}$ and a vector $\phi \in \mathbb{R}^{l_p}$, with $l_p$ and $n$ being positive integers, i.e., }
\begin{equation}
\label{eq:polyset}
   \mathcal{P}(\phi) = \{x \in \mathbb{R}^n: Px \leq \phi\}.
\end{equation}
\end{definition}
Note that $\mathcal{P}(\phi)$ in \eqref{eq:polyset} includes the origin as an interior point iff $\phi > 0$. In the sequel, if $\phi = \textbf{1}_*= [1, \,1, \ldots\,,1\,]^T \in \mathbb{R}^{*}$, the resulting polyhedral set $\mathcal{P}(\textbf{1}_*)$ will be denoted as $\mathcal{P}$.

\begin{definition} 
\label{def:metzler} 
\textit{A matrix $M$ is Metzler type, or essentially non-negative, if $M_{ij}\geq 0, \forall i\neq j$.}
\end{definition}

{\bf \textit{Proof of Lemma 1:}} 
The existence of the Metzler type matrix $H$, the non-negative matrix $H_r$ and the scalar $\gamma>0$ verifying the conditions \eqref{eq:deltainv} are necessary and sufficient algebraic conditions for the robust positive invariance of the set $\mathcal{L}$,  which is the equivalent of imposing the one step admissibility condition $A_{cl}\mathcal{L} \oplus B_{cl} \mathcal{R}(\rho) \subseteq \mathcal{L}$ (see \cite{Lucia2023,castelan1993invariant, blanchini2015set}).

\begin{lemmaBold}\label{lemma:EFL} \textit{(Extended Farkas' Lemma) Consider two polyhedral sets of $\mathbb{R}^n$ defined by $\mathcal{P}_i(\phi_i) =\{x \in \mathbb{R}^n, P_ix \leq \phi_i\}$, for $i=1,2$, with $P_i \in \mathbb{R}^{l_{p_i}\times n}$ and positive vectors $\phi_i \in \mathbb{R}^{l_{p_i}}$. Then, $\mathcal{P}_1 \subseteq \mathcal{P}_2$ if and only if there exists a non-negative matrix $Q \in \mathbb{R}^{l_{p_2} \times l_{p_1}}$ such that}
\begin{equation}
\begin{array}{c}
QP_1 = P_2,\quad  Q\phi_1 \leq \phi_2.
\end{array}
\end{equation}
\end{lemmaBold}

\bibliographystyle{IEEEtran}
\bibliography{Ref}

\end{document}